\documentclass[12pt]{article}

\usepackage{amsmath}
\usepackage{amscd}
\usepackage{amsfonts}
\usepackage{amssymb}
\usepackage{url}
\usepackage{fullpage}

\newtheorem{theorem}{Theorem}%[section]

\newcommand{\diag}{\operatorname{diag}}

\begin{document}

\title{A sneaky proof of the maximum modulus principle}
\author{Orr Moshe Shalit}
\date{}
\maketitle

\begin{abstract}
A proof for the maximum modulus principle (in the unit disc) is presented. This proof is unusual in that it is based on linear algebra. 
\end{abstract}

The goal of this note is to provide a neat proof of the following version of the maximum modulus principle.

\begin{theorem}\label{thm:max}
Let $f$ be a function analytic in a neighborhood of the closed unit disc $\overline{\mathbb{D}} = \{z \in \mathbb{C} : |z|\leq 1\}$. Then 
\[
\max_{z \in \overline{\mathbb{D}}} |f(z)|  = \max_{z \in \partial \mathbb{D}} |f(z)| .
\]
\end{theorem}
(Here and below, $\partial \mathbb{D}$ denotes the unit circle $\partial \mathbb{D} = \{z \in \mathbb{C} : |z| = 1\}$). 

Familiar proofs derive this theorem from the open mapping principle \cite{Ahlfors,GreeneKrantz}, from Cauchy's integral formula via the mean value property for analytic functions \cite{Ash,ChurchillBrown,Kodaira}, or from the maximum principle for subharmonic functions \cite{Evgrafov}. There is also a direct proof which uses the power series representation \cite{Knopp}. 
The proof I will present uses linear algebra, and is motivated by \cite{LevyShalit} and \cite{McCarthyShalit}. 
%Here with page numbers:
%Familiar proofs derive this theorem from the open mapping principle (\cite[p. 170]{GreeneKrantz} or \cite[p. 133]{Ahlfors}), from Cauchy's integral formula via the mean value property for analytic functions (\cite[pp. 117--118]{ChurchillBrown}, \cite[p. 51]{Kodaira} or \cite[p. 42]{Ash}), or from the maximum principle for subharmonic functions (see \cite[p. 224]{Evgrafov}). There is also a direct proof which uses the power series representation \cite[p. 84]{Knopp}. 
%The proof I will present uses linear algebra, and is motivated by the works \cite{LevyShalit} and \cite{McCarthyShalit}.

Before presenting the proof, let me review the main ingredients. For every $x = (x_1, \ldots, x_n) \in \mathbb{C}^n$, we denote 
\[\|x\| = \sqrt{|x_1|^2 + \ldots + |x_n|^2} .
\] 
If $A$ is an $m \times n$ matrix, we define its {\em operator norm} by 
\[
\|A\| = \sup_{\|x\|=1} \|A x\| .
\] 
The only properties of the operator norm that we will require are the following three properties, which are easy consequences of the definition. 
\begin{itemize}
\item[{\bf (a)}] $\|AB\| \leq \|A\| \|B\|$ for all (appropriately sized) matrices $A,B$.
\item[{\bf (b)}] If $D = \diag (d_1, \ldots, d_n)$ then $\|D\| = \max_i |d_i|$. 
\item[{\bf (c)}] If $A$ is unitarily equivalent to $B$, then $\|A\| = \|B\|$. 
\end{itemize}
We will need the following basic fact.

\vskip 5pt

\noindent {\bf Basic Fact 1.} {\em Every unitary matrix is unitarily equivalent to a diagonal matrix $D$, such that the diagonal elements of $D$ all have modulus one.}

\vskip 5pt

It is worth mentioning that Basic Fact 1 does not require any result in complex analysis, and in particular it does not require the fundamental theorem of algebra (see \cite{GG81}, Chapter  III and Section VIII.2). 

The final ingredient we need for the proof is as follows.

\vskip 5pt

\noindent {\bf Basic Fact 2.} {\em If a function is analytic in a neighborhood of $\overline{\mathbb{D}}$, then it is the uniform limit on $\overline{\mathbb{D}}$ of polynomials.}

\vskip 5pt

In fact, a function analytic in a neighborhood of the closed unit disc has a power series representation in some larger disc, and this power series converges uniformly on the closed unit disc (see, e.g., \cite[pp. 79--81]{GreeneKrantz}). 

We are ready to prove Theorem \ref{thm:max}. Note that it suffices to prove the theorem for polynomials. Indeed, suppose the theorem holds for polynomials, let $f$ be as in the theorem and let $\epsilon > 0$. By Basic Fact 2 there is a polynomial $p$ such that $\sup_{\overline{\mathbb{D}}} |f-p|< \epsilon$. Thus
\[
\max_{\overline{\mathbb{D}}} |f| \leq \max_{\overline{\mathbb{D}}} |p| + \max_{\overline{\mathbb{D}}} |f-p| \leq \max_{\partial \mathbb{D}} |p| + \epsilon \leq \max_{\partial \mathbb{D}} |f| + 2 \epsilon ,
\]
whence the theorem follows. We may therefore assume that $f$ appearing in the statement of the theorem is a polynomial. 

Let $n$ be the degree of the polynomial $f$, and let $z$ be any point in $\overline{\mathbb{D}}$. We need to prove that $|f(z)| \leq \max_{\partial \mathbb{D}} |f|$. Put $s = \sqrt{1-|z|^2}$, and define the follwing $(n+1) \times (n+1)$ matrix
\[
U = \begin{pmatrix}
z & & & &  s \\
s & & & &  -\overline{z} \\
& 1 &  & & \\
& & \ddots & & \\
 &  &  & 1 &  
\end{pmatrix} .
\]
%\[
%U = \begin{pmatrix}
%z & 0 & 0 & 0 & \cdots &  0 & s \\
%s & 0 & 0 & 0 &  \cdots & 0 & - \overline{z} \\
% 0 & 1 & 0 & 0 &  & 0 & 0\\
% 0 & 0  & 1 & 0 &    & 0 & 0 \\
% \vdots & \vdots  &   & \ddots & & \vdots & \vdots \\
% \vdots & \vdots  &   &  & 1 & 0 & 0 \\
% 0 & 0  &   & & 0 & 1 & 0
%\end{pmatrix} .
%\]
The empty slots are understood as $0$'s, and the sub-diagonal dots are all $1$'s. One may directly check that $U$ is a unitary matrix. Furthermore, if $P$ denotes the $n+1$ column matrix having $1$ in the $1$st slot and $0$'s elsewhere, and if $P^t$ denotes the transpose of $P$, then a calculation shows that for all $k=1, \ldots, n$,  
\[
z^k = P^t U^k P.  
\]
Since the degree of $f$ is $n$, we obtain that
\[%\label{eq:dil}
f(z) = P^t f(U) P. 
\]
Now, the definition of the operator norm implies that $\|P\| = \|P^t\| = 1$, therefore, using property {\bf (a)} of the operator norm, we find that 
\begin{equation}\label{eq:norm}
|f(z)| \leq \|f(U)\| .
\end{equation}
By Basic Fact 1, $U$ is unitarily equivalent to $\diag(w_1, \ldots, w_{n+1})$, where $|w_i| = 1$ for all $i = 1, \ldots, n+1$. But 
\[
f\Big(\diag(w_1, \ldots, w_{n+1}) \Big) = \diag(f(w_1), \ldots, f(w_{n+1})) ,
\]
thus (using properties {\bf (b)} and {\bf (c)} of the operator norm),
\begin{equation}\label{eq:this}
\|f(U)\| = \|\diag(f(w_1), \ldots, f(w_{n+1}))\| = \max_{1\leq i \leq n+1} |f(w_i)| \leq \max_{\partial \mathbb{D}} |f|. 
\end{equation}
The proof of Theorem \ref{thm:max} is completed by combining equations (\ref{eq:norm}) and (\ref{eq:this}).

\vskip 10pt

\noindent {\bf Concluding Remarks.} 
%\section{Concluding remarks}
\begin{itemize}
\item The proof of Theorem \ref{thm:max} works just as well for functions which are merely continuous on $\overline{\mathbb{D}}$ and analytic in the open disc $\mathbb{D}$. 
%\item The key step in the proof is finding a matrix $U$ which satisfies equation (\ref{eq:dilation}). 
\item There is a simpler proof of the fact that a polynomial cannot achieve its maximum modulus at the {\em center} of a disc. However, to obtain the maximum principle as we stated, one needs an additional ``local to global" argument. The nice thing about the proof given above is that it takes care of the entire disc in one swoop, with no need for an additional argument. 
\item The maximum modulus principle for general bounded domains in $\mathbb{C}$ can be readily deduced from Theorem \ref{thm:max}. However, it would be interesting to find a linear algebra proof along the lines of the above proof that can be applied directly to other domains.
\end{itemize}

\paragraph{Acknowledgments.}  I wish to thank Daniel Reem, for reading a first draft of this paper and offering some very helpful advice.

\bigskip

\noindent\textit{Department of Mathematics\\
Faculty of Natural Sciences\\
Ben-Gurion University of the Negev\\
Be'er Sheva, 84105, Israel \\
oshalit@math.bgu.ac.il} 

\end{document}